\documentclass{amsart}
\usepackage{amsmath,color}
\usepackage{amssymb}
\usepackage{amsfonts}
\usepackage{amsthm}
\usepackage{amscd}
\usepackage{graphicx}
\usepackage{subfigure}
\usepackage{enumerate}

\newtheorem{thm}{Theorem}[section]

\newtheorem{defn}[thm]{Definition}
\renewcommand{\proof}[1][Proof]{\noindent\textsc{#1}. }

\newcommand{\nd}{\noindent}

\begin{document}

\title{ Music By Numbers}

\author{Mihail Cocos\\ Shawn Fowers}

\begin{center}
\begin{figure}
\includegraphics[scale=.2]{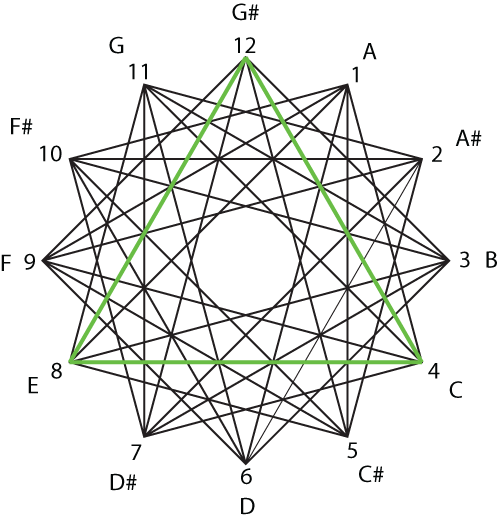} 
\end{figure} 

\end{center} 

\begin{abstract}

In this paper we present a mathematical way of defining musical modes, we derive a formula for the total number of modes and define the musicality of a mode as the total number of harmonic chords whithin the mode. We also give an algorithm for the construction of a duet of melodic lines given a sequence of numbers and a mode. We attach the .mus files of the counterpoints obtained by using the sequence of primes and several musical modes.

\end{abstract}

\maketitle

\section{Musical scales and modes} 
 First we must define some concepts.  In traditional western style of music there are 12 tones in an octave.  The 13th note restarts the cycle and is considered the same note as our first note (simply an octave higher).  Therefore we can assign numerical values to the notes as follows:
 
\begin{center}

 $\begin{array}{ccccc}
 A=1 & A^{\#}=2\\
 B=3 & C=4\\
 C^{\#}=5  & D=6\\
 D^{\#}=7 & E=8 \\
 F=9 & F^{\#}=10\\
 G=11 & G^{\#}=12
 \end{array}$
 \end{center}

 \noindent This gives us a basis to relate musical concepts with mathematics. With this relationship defined, in its most general sense a melodic line simply becomes an extended sequence of numbers.  Most frequently, our choices for which notes are acceptable is dictated by what we call a scale.  
 \begin{defn}
 Definition:  A $k-$scale is an increasing sequence of $k$ integers from one to twelve.  
 \end{defn}

\newpage

\noindent {\bf Examples:} 
\begin{enumerate}
\item  $\mathcal{A}_{min}=(1,3,4,6,8,9,11)$. These numbers will correspond to all the natural tones starting in $A$ and hence we form what is known as the $A$ minor scale.
\item  $\mathcal{A}^{\#}_{ins}=(2,3,7,9,12).$ These numbers correspond to the tones \\ $A^{\#}$,$B$,$D^{\#}$,$F$,$G^{\#}$ to form the $A^{\#}$ In Sen scale.
\end{enumerate}
\begin{defn}
 
Two scales with the same number of tones $S_1=(n_1,n_2, \cdots n_k)$ and $S_2=(m_1,m_2, \cdots  m_k)$ are said to be equivalent iff \[m_i-n_i=m_j-n_j\] for all $1 \leq i,j \leq k.$ 
 \end{defn}
 
\noindent {\bf Note:} 

\begin{enumerate}[(a)]
\item The set of tones whithin a scale define which notes we are allowed to choose when writing a piece of music.  
\item Our definition of a scale will not cover all of the musical scales but will have an equivalent representant.\\ 
\end{enumerate}

\begin{defn}  The set of all equivalent $k-$scales is called a $k-$mode.
 
 \end{defn}
 
Our first goal is to find out how many modes there are. Since the mode of a scale is determined by the spacing between the notes, regardless of which note we begin our scale with, we can assume without loss of generality that each mode begins on $A$ or $1$, and then you choose $k-1$ out of $11$  of the remaining tones in order to produce a mode, therefore we have:
 \begin{thm}
 
   The number of $k-$toned modes is $ C_{11}^{k-1}.$

 \end{thm}

\bigskip                                                                                  
                                                                                  
\nd Crunching the numbers we can calculate that there are a total of $2048$ possible modes.  Taking into account that any givenmusical scale can start on any of our twelve tones, we multiply the number of modes by $12$ to determine that there are $24,576$ different scales!  	Not every one of these scales is musically practical, however.  Not very often do we find a song that is composed entirely using tones $1$ through $7$, for example.  A song composed in this manner would sound rather dissonant to the typical ear.  To mathematically understand why, we should take a look at what is known as a harmonic chord.

\begin{defn}
The interval between two tones is defined as the absolute value of the diference between their numerical values. If the interval is one we call it a semi-tone.

\end{defn}

\begin{defn}  A harmonic interval is an interval of $3, 4, 5, 7, 8,$ or $9$ semitones. 

\end{defn}

\nd Western music has evolved to where these intervals are considered generally more pleasant to the ear than other intervals between notes.  These intervals form the building blocks for major and minor chords.

\begin{defn}  A harmonic subset of a scale $S$ is a subset of tones in $S$ such that the interval between any two of its tones is a harmonic interval.
\end{defn}

\begin{thm}
 The maximum number of tones in a harmonic subset is three.
 \end{thm}

\proof Let $1\leq n_1<n_2< \cdots <n_s \leq 12,$ be the ordered numerical values of the tones in our harmonic subset and assume that $s>3.$ By the very definition of a harmonic subset $|n_i-n_j| \in \{3,4,5,7,8,9\}$ and looking at the first $4$ numerical values of the tones within our set $n_1<n_2<n_3<n_4$ we conclude that $n_4-n_1=9.$ It follows that $n_2-n_1=3, n_3-n_1=3,$ and $n_4-n_3=3,$ which in turn implies
that $n_3-n_1=6 \notin \{3,4,5,7,8,9\},$ which contradicts our hypothesis. $\Box$ \\

\noindent The above theorem allows us to give the following definition.

\begin{defn}

A $3-$toned harmonic subset of a scale is called a chord.

\end{defn}

\nd {\bf Remark:} The reader could easily verify that the only types of harmonic chords $ n_1<n_2<n_3$ are:

\begin{itemize}
\item $n_2-n_1=3$ \& $n_3-n_2=4.$
\item $n_2-n_1=4$ \& $n_3-n_2=3.$
\item $n_2-n_1=3$ \& $n_3-n_2=5.$
\item $n_2-n_1=5$ \& $n_3-n_2=3.$
\item $n_2-n_1=4$ \& $n_3-n_2=4.$
\end{itemize}

\bigskip

\nd Another important fact about two equivalent scales is:

\begin{thm}

Let $S_1=(n_1,n_2, \cdots n_k)$ and $S_2=(m_1,m_2, \cdots  m_k)$ be two equivalent scales (i.e they belong to the same mode), then the number of harmonic chords in $S_1$ equals the number of harmonic chords in $S_2.$

\end{thm}

\proof We may assume, without the loss of generality that $n_i \leq m_i,$ for any $1 \leq i \leq k.$
Let $C_1$ denote the set of harmonic chords in $S_1$ and $C_2$ the set of harmonic chords in $S_2.$
To prove that they have the same cardinality we need to construct a bijection $f$ from $C_1$ to $C_2.$
Let $c=\{n_{i_1},n_{i_2},n_{i_3}\}$ be an arbitrary harmonic chord in $S_1$ and define

\[f(c)= f(\{n_{i_1},n_{i_2},n_{i_3}\})=\{m_{i_1},m_{i_2},m_{i_3}\} .\]

\nd It is easy to see that the intervals between the tones in $f(c)$ are the same as in $c$ and hence $f(c)$ is harmonic and in $S_2.$ It is easy to verify that $f$ has an inverse defined by:

\[ f^{-1}(\{m_{i_1},m_{i_2},m_{i_3}\})=\{n_{i_1},n_{i_2},n_{i_3}\} .\]

$\Box$ \\

\noindent Since all the scales within a mode have the same number of chords we can give the following definition.

\begin{defn}

The musicality of a $k-$mode is the number of harmonic chords of a scale within the mode. 

\end{defn} 
\nd {\bf Note:} This definition only compares modes with the same number of tones.

\bigskip

Now that we know what a harmonic chord is, let's look at the mode containing the first seven semitones.  The only  possible harmonic chord among those tones is $\{1,4,7\}$, or $ \{A, C, D^{\#}\}$ (for the scale begining on $A$).  Contrast this with the commonly used major and minor $7$-toned modes where there are $6$ different harmonic chords and we can see a potential reason why these two modes are used much more than some other $7-$toned modes.

\medskip
\noindent {\bf Open question:}
\begin{center}
{ \em  What are the most musical $k-$toned modes?}
\end{center}
\medskip
	
\nd {\bf Note:}	For more variety, one might consider the interval of six semitones to be included in our possibilities for harmonic intervals.  This interval, when played by itself, isn't the most pleasant of sounds, but using this interval we can add diminished and seventh chords to our types of harmonic chords.  This also increases the maximum number of notes in a harmonic chord from three to four.
To see this take the $1$st, $4$th,$7$th, and $10$th semitones.  The intervals between these notes are three, six, and nine semitones.  Since we are currently considering an interval of six semitones to be harmonic, we have a four toned harmonic chord. \\

\section{The chart of the harmonic chords}
Using these concepts we now describe a geometric algorithm for the construction of all the harmonic chords in any given mode.  The picture on the title page is an example of the picture for the 12-toned mode (chromatic mode!). We present this process as a series of steps.

\bigskip

\begin{itemize}

\item	Begin by assigning numerical values to all 12 tones as we have done on page one of the paper.
\item	Choose the mode to be examined.
\item	Space the notes evenly around a circle along with their given values. Here are examples:

\vspace{0.5 in}
\begin{figure}[h!]
\centering
\mbox{ \subfigure{\includegraphics[scale=0.2]{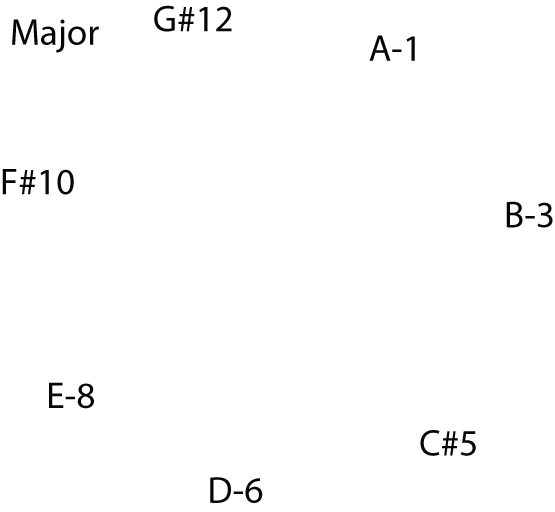}}  \hspace{ 0.5 in}
\subfigure{\includegraphics[scale=0.2]{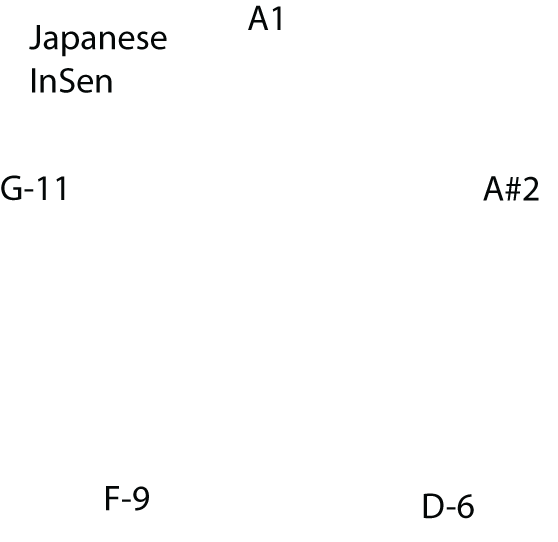} }}
\end{figure}

\newpage

\item 	Draw a line between each of the notes where the difference between the numerical values of the notes is a harmonic interval. 

\vspace{0.5 in}

\begin{figure}[h!]
\centering
\mbox{
\subfigure{\includegraphics[scale=0.2]{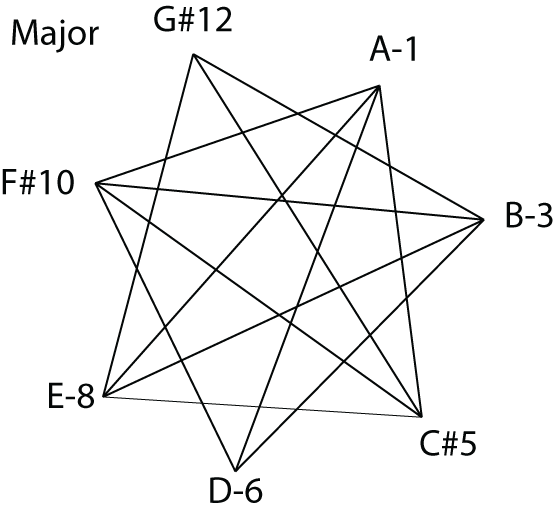}} \hspace{ 0.5 in}  \subfigure{\includegraphics[scale=0.2]{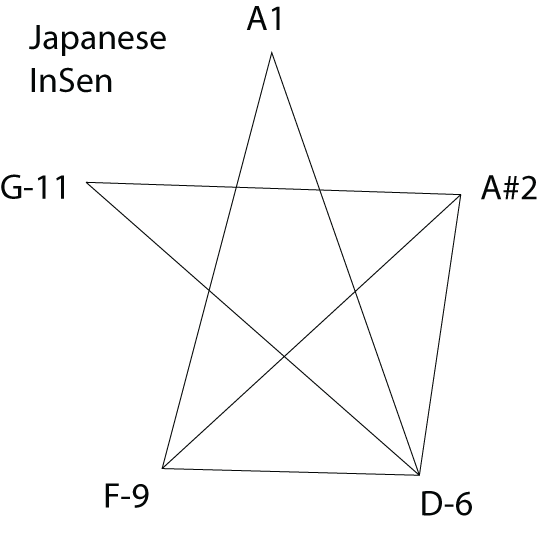}}}
\end{figure}

\item Once we have done this for all of the notes, we examine the picture for any completed triangles.  Each completed triangle represents a harmonic chord.

\vspace{0.5 in}

\begin{figure}[h!]
\centering
\mbox{ \subfigure{\includegraphics[scale=0.2]{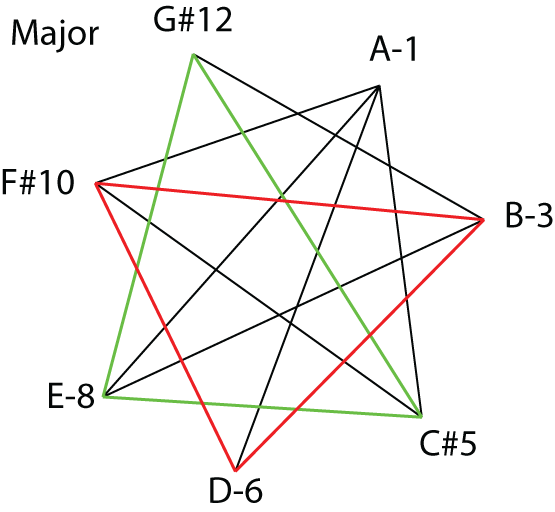}} \hspace{ 0.5 in}
\subfigure{\includegraphics[scale=0.2]{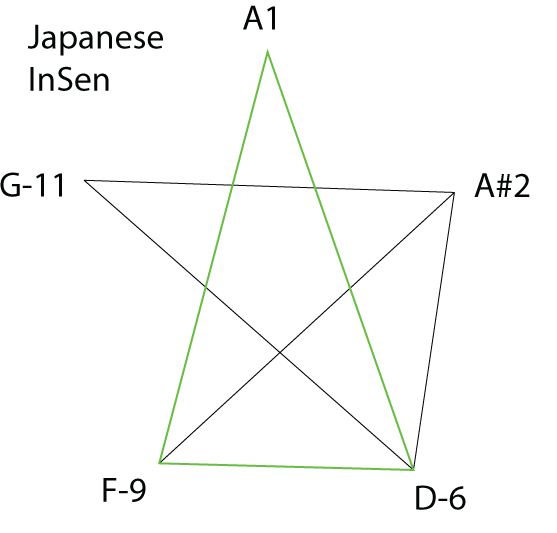}}}
\end{figure}

\end{itemize}

\nd This process will allow us to easily create custom modes and see which harmonic chords are available to us in that mode.

\section{How will a sequence of numbers sound?}

Our final goal in our project is to use our concept of harmonic intervals to compose a song out of a sequence of numbers.  We will use the sequence of prime numbers for this part of the project.  With a list of the first $500$ prime numbers we are ready to begin.
	First we must choose a mode.  Any mode can be selected and the principle behind the technique remains the same.  We will explain our method using the $7$-toned major key, which consists of the notes $A, B, C^{\#}, D, E, F^{\#}$, and $G^{\#}$.  Each of these notes is assigned a  new value from $1$ to $7$ respectively.
	We next compute the prime numbers modulo $7+1=8$ (the modulus will change depending on how many tones are in your mode).  This will give us a value between one and seven, thereby telling us which notes comprise our melody.  Once we have computed this we can enter the notes associated with each prime number modulo $7$ and we have our melody.  
	Next, we take a look at our harmonic intervals  $3, 4, 5, 7, 8,$ and $9$.  We associate each of these with a number $1$ through $6$ the same way we did with the tones of our scale.  Computing the prime numbers modulo $6+1=7$ this time we can now associate each prime number with a harmonic interval.  Once we have found our intervals, we take a look at the melody note that is associated with each prime number and we subtract the value of the harmonic interval from the value that we gave each note of the $12$-toned octave on page $1$.  This new number is the number of our potential harmony note.  If the result of this calculation is negative, simply extend the $12-$toned scale downward into negative values (the $12$ tones repeat themselves as high and as low as you wish to go).  
	We now examine the harmony notes.  If it is included in the notes of the mode that we chose then we add this to our harmony line.  If the note is not included in our mode then we leave that beat as a rest.   We now have a basic duet.  It may sound strange in some modes while it may sound better in other modes, but this technique will work with any sequence of numbers in any mode.  Try different ones until you find something you like.
\nd We may use this method to assign rhythmic values to our notes as well.  To do this we choose which values of notes we wish to include in our song (e.g.  quarter note, eighth note, half note) and assign a number from $1$ to $n$ to each duration where $n$ is the number of different rhythmic values we have chosen.  Computing the prime numbers modulo $n+1$ we are returned a value representative of the rhythmic value of the corresponding notes.  Using this method may work in some cases (such as the blues number we have included) where the rhythms are often syncopated, however, in other cases this method does not lend itself well to a metered song.  Another attempt (not included) proved to be terribly erratic and difficult for the ear to follow.  Perhaps we can confine our choices for rhythmic values to ones that will give us a greater chance to finish out measures evenly, or perhaps you can think of another way that would be more successful for assigning rhythmic values using our number sequence?

\end{document}